\documentclass[11pt]{amsart}
\usepackage{amssymb,latexsym,eucal,mathrsfs,array}

\usepackage[curve,v2]{xypic}

\advance\textwidth by 1.2in
\advance\oddsidemargin by -.6in
\advance\evensidemargin by -.6in

\newcommand{\tensor}{\otimes}

\newcommand{\cliff}{\operatorname{Cliff}}
\newcommand{\picard}{\operatorname{Pic}}
\newcommand{\ext}{\operatorname{Ext}}

\newcommand{\wt}[1]{\widetilde{#1}}
\newcommand{\bt}{\boxtimes}
\newcommand{\wtp}{\widetilde{\P^n}}
\newcommand{\wtv}{\widetilde{\varphi }}
\newcommand{\wts}{\widetilde{Sec X}}
\newcommand{\wtvp}[1]{\widetilde{\varphi_{#1}}^+}
\newcommand{\wtvn}[1]{\widetilde{\varphi_{#1}}^-}

\newcommand{\posmap}[1]{\wtvp{#1}:M_{#1}\rightarrow \P^{s_{#1}}}
\newcommand{\negmap}[2]{\wtvn{#2}:M_{#1}\rightarrow \P^{s_{#2}}}
\newcommand{\PP}[1]{\P^{s_{#1}}}

\newcommand{\blow}[2]{{\rm Bl}_{#2}{(#1})}

\newcommand{\ses}[3]{0\rightarrow#1\rightarrow#2
   \rightarrow#3\rightarrow0}

\newcommand{\U}{{\mathcal U}}

\renewcommand{\H}{{\mathcal H}}
\newcommand{\F}{{\mathscr F}}
\newcommand{\D}{{\mathcal D}}
\newcommand{\E}{{\mathscr E}}

\newcommand{\I}{{\mathscr I}}
\renewcommand{\L}{{\mathscr L}}
\renewcommand{\O}{{\mathcal O}}
\renewcommand{\P}{{\mathbb{P}}}
\newcommand{\C}{{\mathbb{C}}}

\newcommand{\Z}{{\mathbb{Z}}}
\newcommand{\Q}{{\mathbb{Q}}}

\newcommand{\G}{{\mathbb{G}}}

\newcommand{\B}{{\mathcal{B}}}

\begingroup
\newtheorem{thm}{Theorem}[section]   
\newtheorem{lemma}[thm]{Lemma}         
\newtheorem{prop}[thm]{Proposition}  
\newtheorem{defn}[thm]{Definition}   
\newtheorem{conj}[thm]{Conjecture}        
   
\endgroup

\newenvironment{rem}[2]{\refstepcounter{thm} \label{#2} 
\par \medskip \noindent {\bf #1 \thethm }}{\par \medskip}

\newcommand{\gmodh}[2]{\raisebox{.1cm}{\ensuremath{#1}}/\raisebox{-.1cm}{\ensuremath{#2}}} 

\newenvironment{pic}[3]{
\input epsf
\begin{figure}[htb] 
  \begin{center}
     \leavevmode
      \hbox{%
      \epsfxsize= #2 in
      \epsfysize= #3 in
      \epsffile{#1}}
   \end{center}
\end{figure}
}

\begin{document}

\pagenumbering{arabic}

\title[]{Secant Varieties and Birational Geometry}

\author[]{Peter Vermeire}

\address{Department of Mathematics, Oklahoma State
University, Stillwater OK 74078}

\email{petvermi@math.okstate.edu}

\date{\today}

\begin{abstract}
We show how to use information about the equations defining
secant varieties to smooth projective varieties in order
to construct a natural collection of birational
transformations.  These were
first constructed as flips in the case of curves by M. Thaddeus via
Geometric Invariant Theory, and the first flip in the sequence 
was constructed by the author for varieties of arbitrary dimension 
in an earlier paper.  We expose the finer structure
of a second flip; again for varieties of arbitrary dimension.
We also prove a result on the cubic generation of the secant 
variety and give some conjectures
on the behavior of equations defining the higher secant varieties.
\end{abstract}

\maketitle

\section{Introduction}
In this paper we continue the geometric construction of a sequence of
flips associated to an embedded projective variety begun in 
\cite{verm-flip1}.  We give hypotheses under which this sequence
of flips exists, and state some conjectures on how positive a line
bundle on a curve must be to satisfy these hypotheses.  These conjectures
deal with the degrees of forms defining various secant varieties
to curves and seem interesting outside of the context of the flip 
construction.

As motivation, we have the work of A. Bertram and M. Thaddeus.
In \cite{thad-flip} this sequence of flips is constructed in the
case of smooth curves via GIT, in the context of the moduli space
of rank two vector bundles on a smooth curve.  An understanding
of this as a sequence of log flips is given in \cite{bertram3}, and
further examples of sequences of flips of this type, again 
constructed via GIT, are given in \cite{thad-toric},\cite{thad-coll}.
Our construction, however, does not use the tools of Geometric
Invariant Theory and is closer in spirit to \cite{bertram1},\cite{bertram2}.

In Section~\ref{review}, we review the constructions in \cite{bertram1}
and \cite{thad-flip} and describe the relevant results from \cite{verm-flip1}.
In Section~\ref{genofsecants} we discuss the generation of $SecX$
by cubics.  In particular, we show (Theorem~\ref{settheo}) that
large embeddings of varieties have secant varieties that are
at least set theoretically defined by cubics.  We also offer
some general conjectures and suggestions in this direction for
the generation of higher secant varieties.

The construction of the new flips is somewhat more involved than that
of the first in \cite{verm-flip1}.  We give a general construction
of a sequence of birational transformations in
Section~\ref{genflipconst}, and we describe in detail the 
second flip in Section~\ref{secondflip}.

We mention that some of the consequences of these constructions and
this point of view are worked out in \cite{verm-vanish}.

{\bf Notation:} We will decorate a projective variety $X$ as follows: 
$X^d$ is the $d^{th}$ cartesian product of $X$; $S^dX$ is
$Sym^dX=\gmodh{X^d}{S_d}$, the $d^{th}$ symmetric product of $X$;
and $\H^dX$ is $Hilb^d(X)$, the Hilbert Scheme of zero dimensional
subschemes of $X$ of length $d$.  Recall (Cf. \cite{gott}) that if $X$ is
a smooth projective variety then $\H^dX$ is also projective, and is smooth
if and only if either $dim~X\leq 2$ or $d\leq 3$.

Write $Sec_{k}^{\ell}X$ for the (complete) variety of $k$-secant
$\ell$-planes to $X$.  As this notation can become cluttered,
we simply write $Sec^{\ell}X$ for $Sec^{\ell}_{\ell +1}X$
and $SecX$ for $Sec^1_2X$.  Note also the convention $Sec^0X=X$.
If $V$ is a $k$-vector space, we denote by $\P(V)$ the space of
1-dimensional quotients of $V$.  Unless otherwise stated, we work
throughout over the field $k=\C$ of complex numbers.
We use the terms locally free sheaf (resp. invertible sheaf) and vector 
bundle (resp. line bundle) interchangeably.
Recall that a line bundle $\L$ on $X$ is {\em nef} if $\L.C\geq 0$ for every 
irreducible curve $C\subset X$.  A line bundle $\L$ is {\em big} if 
$\L^{\tensor n}$ induces a birational map for all $n\gg 0$.

{\bf Acknowledgments:} I would like to thank Aaron Bertram, Sheldon
Katz, Zhenbo Qin, and Jonathan Wahl for their helpful conversations
and communications.

\section{Overview of Stable Pairs and the Geometry of SecX}\label{review}

Fix a line bundle $\Lambda$ on a fixed smooth curve $X$, and 
denote by $M(2,\Lambda)$ the moduli space of semi-stable rank two vector
bundles $E$ with $\wedge^2E=\Lambda$.
There is a natural rational map, the {\it Serre Correspondence}
$$\Phi:\P(\Gamma(X,K_X\tensor\Lambda)^*)\dashrightarrow M(2,\Lambda)$$
given by the duality $\ext^1(\Lambda,\O)\cong
H^1(X,\Lambda^{-1})\cong H^0(X,K_X\tensor \Lambda)^*$, taking
an extension class $\ses{\O}{E}{\Lambda}$ to $E$.  One has
an embedding $X\hookrightarrow\P(\Gamma(X,K_X\tensor\Lambda)^*)$
(at least in the case $d=c_1(\Lambda)\geq 3$) and 
$\Phi$, defined only for semi-stable $E$, is a morphism off $Sec^kX$ where 
$k=\left[ \frac{d-1}{2}\right] $ \cite{bertram2}. 
This map is resolved in \cite{bertram1} by first blowing up along $X$,
then along the proper transform of $SecX$, then along the
transform of $Sec^2X$ and so on until we have a morphism to $M(2,\Lambda)$.

A different approach is taken in \cite{thad-flip}.  There, for a
fixed smooth curve $X$ of genus at least $2$ and a fixed line bundle $\Lambda$,
the moduli problem of semi-stable pairs $(E,s)$ consisting of a rank 
two bundle $E$ with $\wedge^2E=\Lambda$, and a section 
$s\in \Gamma(X,E)-\{ 0\}$, 
is considered.  This, in turn, is interpreted as a GIT problem, and
by varying the linearization of the group action, a collection
of (smooth) moduli spaces $M_1,M_2, \ldots,M_k$ ($k$ as above) 
is constructed.  As stability is an open condition, these spaces
are birational.  In fact, they are isomorphic in codimension one, and
may be linked via a diagram
\begin{center}
{\begin{minipage}{1.5in}
\diagram
 & \wt{M_2} \dlto \drto & & \wt{M_3} \dlto \drto & & \wt{M_k} \dlto \drto & \\
 M_1  & & M_2 & & \cdots & & M_k
\enddiagram
\end{minipage}}
\end{center}
where there is a morphism $M_k\rightarrow M(2,\Lambda)$.
The relevant observations are first that this is a diagram of flips (in fact
it is shown in \cite{bertram3} that it is a sequence of log flips)
where the ample cone of each $M_i$ is known.
Second, $M_1$ is the blow up of $\P(\Gamma(X,K_X\tensor\Lambda)^*)$
along $X$, $\wt{M_2}$ is the blow up of $M_1$ along
the proper transform of the secant variety, and all of the
flips can be seen as blowing up and down various higher secant varieties.
Finally, the $M_i$ are isomorphic off loci which are 
projective bundles over appropriate symmetric products
of $X$.

Our approach is as follows:  The sequence of flips
in Thaddeus' construction can be realized as a sequence of 
geometric constructions depending
only on the embedding of $X\subset \P^n$.  An
advantage of this approach is that the smooth curve $X$ can
be replaced by any smooth variety.  Even in the curve case, our 
approach applies to situations where Thaddeus' construction does not 
hold (e.g. for canonical curves with $\cliff X>2$).
In \cite{verm-flip1}, we show how to construct the first flip
using only information about the syzygies among the equations
defining the variety $X\subset \P^n$.  We summarize this construction
here.

\begin{defn}\label{kd}
Let $X$ be a subscheme of $\P^n$.  The pair $(X,F_i)$ 
{\bf \boldmath satisfies condition 
$(K_d)$} if $X$ is scheme theoretically cut out by forms 
$F_0, \ldots,F_s$ of degree $d$ such that the trivial (or Koszul)
relations among the $F_i$ are generated by linear syzygies.
\end{defn}
We say $(X,V)$ satisfies $(K_d)$ for $V\subseteq H^0(\P^n,\O(d))$ if
$V$ is spanned by forms $F_i$ satisfying the above condition.  We say 
simply $X$ satisfies $(K_d)$ if there exists a set $\{ F_i\} $ 
such that $(X,F_i)$ satisfies $(K_d)$, and if the discussion depends
only on the existence of such a set, not on the choice of a particular set.

As $(K_2)$ is a weakening of Green's property $(N_2)$\cite{mgreen},
examples of varieties satisfying $(K_2)$ include smooth curves embedded
by complete linear systems of degree at least $2g+3$, canonical curves
with $\cliff X\geq 3$, and sufficiently large embeddings of arbitrary
projective varieties.

To any projective variety $X\subset \PP{0}$ defined (as a scheme) by forms 
$F_0, \ldots,F_{s_1}$ of degree $d$, there is an associated rational map
$\varphi:\PP{0}\dashrightarrow\PP{1}$ defined off the common zero locus
of the $F_i$, i.e. off $X$.  This map may be resolved to a morphism
$\wtv:\wt{\PP{0}}\rightarrow\PP{1}$ by blowing up $\P^n$ along $X$, or
equivalently by projecting from the closure of the graph
$\overline{\Gamma_{\varphi}}\subset \PP{0}\times\PP{1}$.  We have
the following results on the structure of $\wtv$:
\begin{thm}\label{embed}\cite[2.4-2.10]{verm-flip1}
Let $(X,F_i)$ be a pair that satisfies $(K_d)$.  Then:
\begin{enumerate} 
\item $\varphi :\PP{0}\setminus X\rightarrow \PP{1}$ is an embedding 
off of $Sec^1_dX$, the variety of $d$-secant lines.
\item The projection of a positive dimensional fiber of $\wtv$ to 
$\PP{0}$ is either contained in a linear subspace of $X$ or is a linear space
intersecting $X$ in a $d$-tic hypersurface.
\end{enumerate}
If, furthermore, $X$ does not contain a line then $\wtv$ is an 
embedding off the proper transform of $Sec^1_dX$.
{\nopagebreak \hfill $\Box$ \par \medskip}
\end{thm}

\begin{thm}\label{getabundle}\cite[3.8]{verm-flip1}
Let $(X,V)$ satisfy $(K_2)$ and assume $X\subset \PP{0}$ is smooth,
irreducible, contains no lines and contains no quadrics.  Then:
\begin{enumerate}
\item The image of $\wts=\wt{Sec^1_2X}$ under $\wtv$ is $\H^2X$.
\item $\E=\wtv_*(\O_{\wts}(H))$ is a rank two vector bundle
on $\H^2X$, where $\O_{\wt{\PP{0}}}(H)$ is the proper transform of the
hyperplane section on $\PP{0}$.
\item $\wtv:\wts\rightarrow\H^2X$ is the $\P^1$-bundle 
$\P_{\H^2X}\E\rightarrow \H^2X$.
{\nopagebreak \hfill $\Box$ \par \medskip}
\end{enumerate}
\end{thm}

This implies $\wts$, and hence $\wt{M_2}=\blow{\wt{\PP{0}}}{\wts}$, are smooth.
To complete the flip, we construct a base point free linear
system on $\wt{M_2}$, and take $M_2$ to be the image of the associated
morphism.  Denoting $\wts=\P(\E)$, the sheaf 
$\F=\wtv_*(N^*_{\P(\E)/\wt{\PP{0}}}\tensor \O_{\P(\E)}(-1))$ is 
locally free of rank $n-2\dim X-1$ on $\H^2X$.
Write $\P(\F)=\P_{\H^2X}(\F)$ and rename $\wtv$ as $\wtvp{1}$:
\renewcommand{\wtp}{\widetilde{\PP{0}}}
\begin{thm}\label{flip}\cite[4.13]{verm-flip1}
Let $(X,V)$ satisfy $(K_2)$ and assume $X\subset \PP{0}$ is smooth,
irreducible, contains no lines and contains no plane quadrics.  Then
there is a flip as pictured below with:
\begin{enumerate}
\item $\wtp$, $\wt{M_2}$, and $M_2$ smooth
\item $\wtp\setminus \P(\E)\cong M_2\setminus \P(\F)$, hence 
if $codim(\P(\E),\wtp)\geq 2$ then $\picard\wtp\cong \picard M_2$
\item $h_1$ is the blow up of $M_2$ along $\P(\F)$
\item $\pi$ is the blow up of $\wtp$ along $\P(\E)$
\item $\wtvn{1}$, induced by $\O_{M_2}(2H-E)$, is an embedding off
of $\P(\F)$, and the restriction of $\wtvn{1}$ is
the projection $\P(\F)\rightarrow \H^2X$ 
\item $\wtvp{1}$, induced by $\O_{\wtp}(2H-E)$,
is an embedding off of $\P(\E)$, and the restriction of
$\wtvp{1}$ is the projection $\P(\E)\rightarrow \H^2X$
\end{enumerate}
\end{thm}

\begin{center}
{\begin{minipage}{1.5in}
\diagram
 & E_2 \dlto_{\pi} \drto^{h_1} & & & \wt{M_2} \dlto_{\pi}
\drto^{h_1} & \\
  \P(\E) \drto_{\wtvp{1}} & & \P(\F)
\dlto^{\wtvn{1}} &  \wtp \drto^{\wtvp{1}} \dto & & M_2
\dlto^{\wtvn{1}} \\
 & \H^2X & & \P^{s_0} \rdashed|>\tip_{\varphi_1}  & \P^{s_1} &
\enddiagram
\end{minipage}}
\end{center}

To continue this process following Thaddeus, we need to construct 
a birational morphism $\posmap{2}$ which contracts the transforms 
of $3$-secant $2$-planes to points, and is an embedding off their union.  
The natural candidate is the map induced by the
linear system $\O_{M_2}(3H-2E)$.  We discuss two different 
reasons for this choice that will guide the construction of the 
entire sequence of flips.  Section~\ref{genofsecants}
addresses the question of when this system is globally generated.
Note that we abuse notation throughout and identify
line bundles via the isomorphism $\picard \wtp \cong
\picard M_k$.

The first reason is quite naive: Just as quadrics collapse
secant lines because their restriction to such a line is 
a quadric hypersurface, so too do cubics vanishing twice
on a variety collapse every $3$-secant $\P^2$ because they vanish
on a cubic hypersurface in such a plane.  Similarly, to collapse the
transform of each $k+1$-secant $\P^k$ 
via a morphism $\posmap{k}$, the natural system is $\O_{M_k}((k+1)H-kE)$.

Another reason is found by studying the
ample cones of the $M_i$.  Note that the ample cone on 
$\wtp (=M_1)$ is bounded by the line bundles $\O_{\wtp}(H)$ and
$\O_{\wtp}(2H-E)$.  Both of these bundles are globally generated, and
by Theorems~\ref{embed} and \ref{getabundle}, they each
give birational morphisms whose exceptional loci are 
projective bundles over Hilbert schemes of points of $X$ 
($\H^1X\cong X$ and $\H^2X$ respectively).

On $M_2$, the ample cone is bounded on one side by 
$\O_{M_2}(2H-E)$.  This gives the map $\negmap{2}{1}$ mentioned
in Theorem~\ref{flip}; in particular it is globally
generated, the induced morphism is birational, and its 
exceptional locus is a projective bundle over $\H^2X$.
On the other side, the ample cone contains a line bundle
of the form $\O_{M_2}((2m-1)H-mE)$ (\cite[4.9]{verm-flip1}).  In fact,
if $X$ is a smooth curve embedded by a line bundle of degree at
least $2g+5$, it is shown in \cite{thad-flip} that the case $m=2$ 
suffices, i.e. that the ample cone is bounded by $\O_{M_2}(2H-E)$
and $\O_{M_2}(3H-2E)$.  Therefore, it is natural to look for conditions
under which $\O_{M_2}(3H-2E)$ is globally generated.  Thaddeus
further shows that under similar positivity conditions, the
ample cone of $M_k$ is bounded by $\O_{M_k}(kH-(k-1)E)$ and
$\O_{M_k}((k+1)H-kE)$. 

Noting the fact that $h_1^*\O_{M_2}(3H-2E)=\O_{\wt{M_2}}(3H-2E_1-E_2)$, 
it is not difficult to see (using Zariski's Main Theorem) that
this system will be globally generated if $SecX\subset \PP{0}$ is
scheme theoretically defined by cubics, because a cubic vanishing
twice on a variety must also vanish on its secant variety.
Unfortunately, there are no general theorems on the cubic generation
of secant varieties analogous to quadric generation
of varieties.  We address this question in the next section.

\section{Cubic Generation of Secant Varieties}\label{genofsecants}

\begin{rem}{Example}{somecubics}
Some examples of varieties whose secant varieties are {\em ideal}
theoretically defined by cubics include:
\begin{enumerate}
\item $X$ is any Veronese embedding of $\P^n$ \cite{kanev} 
\item $X$ is the Pl\"ucker embedding of the Grassmannian $\G(1,n)$ for any $n$ 
\cite[9.20]{harris}.
\item $X$ is the Segre embedding of $\P^n\times\P^m$ \cite[9.2]{harris}.
{\nopagebreak \hfill $\Box$ \par \medskip}
\end{enumerate}
\end{rem}
We prove a general result:

\begin{thm}\label{settheo}
Let $X\subset \PP{0}$ satisfy condition
$(K_2)$.  Then $Sec(v_d(X))$ is {\bf set} theoretically defined by cubics
for $d\geq 2$.
\end{thm}

\begin{proof}
We begin with the case $d=2$, the higher embeddings being more
elementary.

Let $Y=v_2(X)$, $V=v_2(\PP{0})\subset \P^N$, and $H$ the linear subspace of
$\P^N$ defined by the hyperplanes corresponding to all the quadrics in $\PP{0}$
vanishing on $X$.  Then $Y=V\cap H$ as schemes and we show, noting
that $SecV$ is ideal theoretically defined by cubics, that 
$SecY=SecV\cap H$ as sets.

Note that the map $\varphi_1:\PP{0}\dashrightarrow \P^{s_1}$ 
can be viewed as the
composition of the embedding $v_2:\PP{0}\hookrightarrow \P^N$ with
the projection from $H$, $\P^N\dashrightarrow \P^{s_1}$. 

Let $p\in SecV\cap H$.
If $p\in V$, then $p\in Y=V\cap H$ hence $p\in SecY$.  

Otherwise, any secant line $L$ to $V$ through $p$ intersects $V$
in a length two subscheme $Z$.  $Z$ considered in $\PP{0}$ determines a 
unique line in $\PP{0}$
whose image in $\P^N$ is a plane quadric $Q\subset V$ spanning a 
plane $M$.
If $H\cap Q=Z'\subset Y$ is non-empty then 
$Z'\cup\{ p\} \subset H\cap M$, hence either $H$ intersects $M$
in a line $L'$ through $p$ or $M\subset H$.  In the first case $L'$
is a secant line to $Y$, in the second $Q\subset Y$.  In either 
situation $p\in SecY$.

\begin{pic}{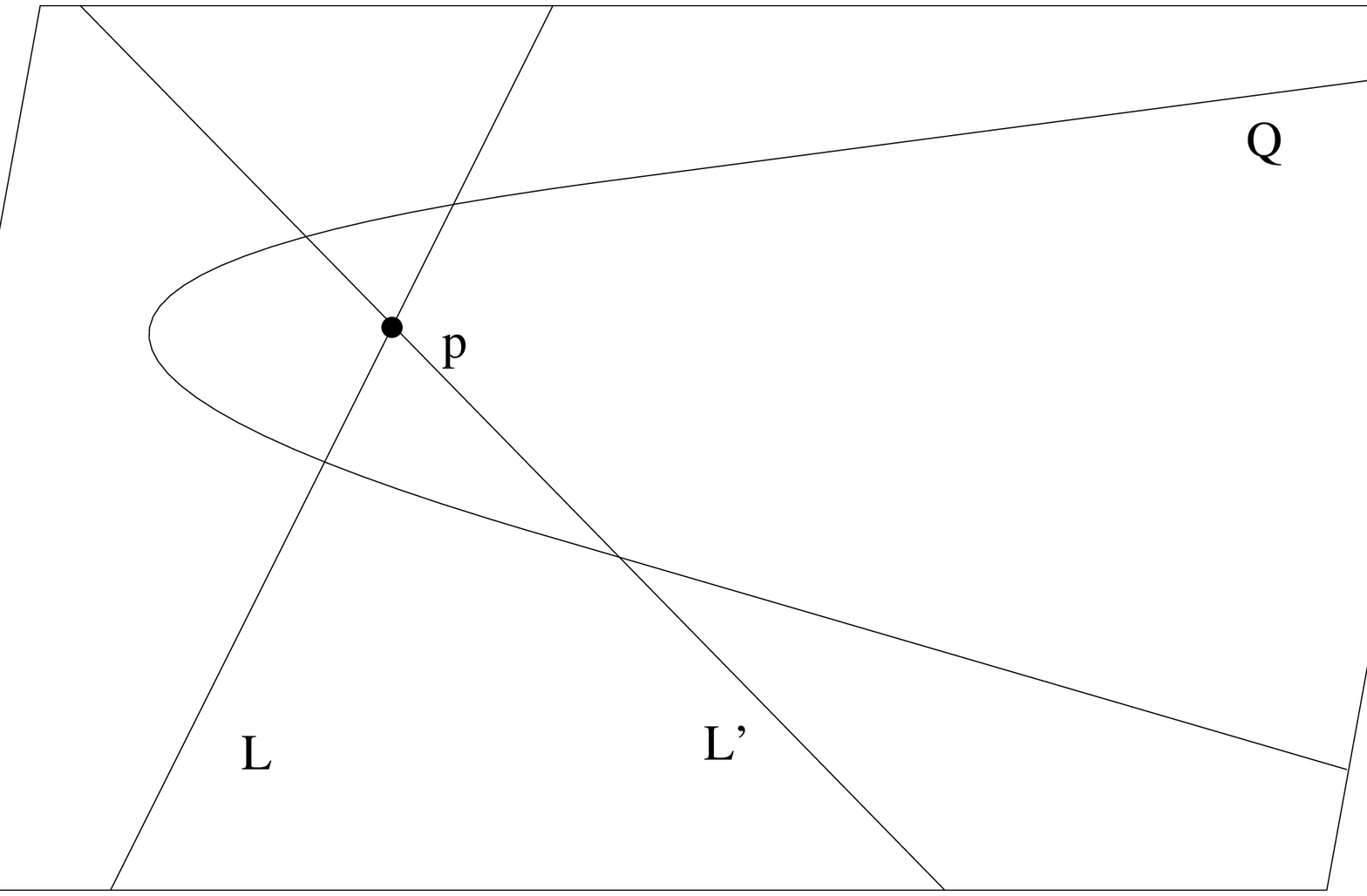}{3}{1.5}
\end{pic}

All that remains is the case $H\cap M=\{ p\} $ and $H\cap Q$
is empty.  However in this case the line $L$, and hence the scheme
$Z=L\cap Q$ is collapsed to a point by the projection.  As the 
rational map $\varphi$ is an
embedding off $SecX$, this implies $Z$ lies on the image of a secant
line to $X\subset \PP{0}$.  As a length two subscheme of $\PP{0}$
determines a unique line, $Q$ must be the image of a secant line
to $X\subset \PP{0}$ contradicting the assumption that $H\cap Q$
is empty.

For $d>2$, note that the projection from $H$ 
is an embedding off $V\cap H$ (this can be derived
directly from Theorem~\ref{embed} or see \cite[3.3.1]{verm-thesis}).
Therefore, if $H$ intersects a secant line, the line lies in $H$, hence
is a secant line to $Y$.
\end{proof}

\begin{rem}{Example}{cubics} As Green's $(N_2)$ 
implies $(K_2)$, this shows that the secant
varieties to the following varieties are set theoretically defined by
cubics:
\begin{enumerate}
\item $X$ a smooth curve embedded by a line bundle of degree $4g+6+2r$,
$r\geq 0$.
\item  $X$ a smooth curve with $\cliff X>2$, embedded by 
$K_X^{\tensor r}$, $r\geq 2$.
\item $X$ a smooth variety embedded by 
$\left( K_X\tensor L^{\tensor (\dim X+3+\alpha)}\right) ^{\tensor r}$, 
$\alpha\geq 0, r\geq 2$, $L$ very ample.
\item $X$ a smooth variety embedded by $L^{\tensor 2r}$ for all $r\gg 0$,
$L$ ample.
{\nopagebreak \hfill $\Box$ \par \medskip}
\end{enumerate}
\end{rem}

\begin{rem}{Remark}{almostk3}
Notice that in the case $d=2$ of Proposition~\ref{settheo}, 
the cubics that at least set theoretically 
define the secant variety satisfy $(K_3)$.  This is because:
\begin{enumerate}
\item The ideal of the secant variety of $v_2(\PP{0})$ is generated 
by cubics, and the module of syzygies is generated by linear 
relations \cite[3.19]{jpw}.  Hence $Sec(v_2(\PP{0}))$ satisfies
$(K_3)$.
\item It is clear from the definition that if $X\subset\P^n$
satisfies $(K_d)$, then any linear section does as well.
{\nopagebreak \hfill $\Box$ \par \medskip}
\end{enumerate}
\end{rem}

\begin{rem}{Example}{quadhyp}
If $X\subset\P^n$ is a smooth quadric hypersurface, then
$v_2(X)$ is given by the intersection of $v_2(\P^n)$ with
a hyperplane $H$.  Furthermore, the intersection of $Sec(v_2(\P^n))$
with $H$ is a scheme $S$ with $S_{red}\cong Sec(v_2(X))$.  
Therefore, a general smooth quadric hypersurface has
$Sec(v_2(X))\cong Sec(v_2(\P^n))\cap H$ as schemes, hence
$Sec(v_2(X))$ satisfies $(K_3)$.
{\nopagebreak \hfill $\Box$ \par \medskip}
\end{rem}

We record here a related conjecture of Eisenbud, Koh, and Stillman as well as
a partial answer proven by M.S. Ravi:

\begin{conj}\cite{eks}
Let $L$ be a very ample line bundle that embeds a smooth curve $X$.  For each
$k$ there is a bound on the degree of $L$ such that $Sec^kX$ is 
ideal theoretically defined by the $(k+2)\times (k+2)$ minors of a matrix of
linear forms.
\end{conj}

\begin{thm}\cite{ravi}
If $deg~L\geq 4g+2k+3$, then $Sec^kX$ is
{\bf set} theoretically defined by the $(k+2)\times (k+2)$ minors of a 
matrix of linear forms.
\end{thm}

These statements provide enough evidence to make the following basic:
\begin{conj}
Let $L$ be an ample line bundle on a smooth variety $X$, $k\geq 1$ fixed.
Then for all $n\gg 0$, $L^n$ embeds $X$ so that $Sec^kX$ is
ideal theoretically defined by forms of degree $k+2$, and furthermore
satisfies condition $(K_{k+2})$.
\end{conj}

\begin{rem}{Remark}{5secant3plane}
If $X$ is a curve with a 
$5$-secant $3$-plane, then any cubic
vanishing on $SecX$ must vanish on that $3$-plane.  Hence
$SecX$ cannot be set theoretically defined by cubics.
This should be compared to the fact that if $X$ has a trisecant line, then
$X$ cannot be defined by quadrics.  In particular, this shows that Green's
condition $(N_2)$ is not even sufficient to guarantee that their exists
a cubic vanishing on $SecX$. For example, if $X$ is an elliptic
curve embedded in $\P^4$ by a line bundle of degree $5$, then 
$SecX$ is a quintic hypersurface.
Therefore, any uniform bound on the degree of a linear system
that would guarantee $SecX$ is even set theoretically defined by cubics
must be at least $2g+4$.
{\nopagebreak \hfill $\Box$ \par \medskip}
\end{rem}

We can use earlier work to give a more geometric necessary condition
for $SecX$ to be defined as a scheme by cubics.  Specifically, in 
\cite[3.7]{verm-flip1} it is shown that the intersection of
$\wts$ with the exceptional divisor $E$ of the blow up of $\PP{0}$
along $X$ is isomorphic to $\blow{X\times X}{\Delta}$.  This implies 
that if $\pi:\wts \rightarrow SecX$ is the blow up along $X$, then $\pi^{-1}(p)
\cong \blow{X}{p}, p\in X$.  In fact, it is easy to verify 
that if $X$ is embedded by a line bundle $L$, then 
$\pi^{-1}(p) \cong \blow{X}{p}\subset \P\Gamma(X,L\tensor \I_p^2)$
where $\P\Gamma(X,L\tensor \I_p^2)$ is identified with the fiber over 
$p$ of the 
projectivized conormal bundle of $X\subset \PP{0}$.  Now, if $SecX$
is defined as a scheme by cubics, then the base scheme of 
$\O_{\wtp}(3H-2E)$ is precisely $\wts$.  The restriction of this series
to $\P\Gamma(X,L\tensor \I_p^2)$ is thus a system of quadrics whose 
base scheme
is $\blow{X}{p}$.  In other words, if $X$ is a smooth variety embedded by a 
line bundle
$L$ that satisfies $(K_2)$ and if $SecX$ is scheme theoretically defined
by cubics, then for every $p\in X$ the line bundle 
$L\tensor \O(-2E_p)$ is very ample on $\blow{X}{p}$
and $\blow{X}{p}\subset \P\Gamma(\blow{X}{p},L\tensor \O(-2E_p))$ is 
scheme theoretically 
defined by quadrics.

In the case $X$ is a curve, this implies that a uniform 
bound on $deg~L$ that would imply $SecX$ is defined by cubics must be at
least $2g+4$, the same bound encountered in Remark~\ref{5secant3plane}.
The construction in \cite{bertram1} shows similarly 
that any uniform bound
that would imply $Sec^kX$ is defined by $(k+2)$-tics must be
at least $2g+2+2k$.  We combine these observations with the degree
bounds encountered in the constructions of \cite{thad-flip} 
and \cite{bertram1} to form the following:

\begin{conj}\label{degreeboundsconj}
Let $X$ be a smooth curve embedded by a line bundle $L$.  If 
$deg(L)\geq 2g+2k$ then $Sec^{k-1}X$ is 
defined as a scheme by forms of degree $k+1$.  If $deg(L)\geq 2g+2k+1$ then 
$Sec^{k-1}X$ satisfies condition $(K_{k+1})$.
{\nopagebreak \hfill $\Box$ \par \medskip}
\end{conj}

\section{The General Birational Construction}\label{genflipconst}

Suppose that $X$ satisfies $(K_2)$, is smooth, and contains no lines
and no plane quadrics.  Suppose further that $SecX$ is scheme theoretically
defined by cubics $C_0, \ldots,C_{s_2}$, and that $SecX$ satisfies $(K_3)$.  
Under these hypotheses, we construct a second flip as follows:  We know
that $\O_{M_2}(3H-2E)$ is globally generated by the discussion above; hence
this induces a morphism $\posmap{2}$ which agrees with the map given 
by the cubics $\varphi_2:\PP{0}\dashrightarrow \P^{s_2}$ on the locus
where $M_2$ and $\PP{0}$ are isomorphic.  By Theorem~\ref{embed}, 
$\wtvp{2}$ is a birational morphism.
We wish first to identify the exceptional locus of $\wtvp{2}$.  It is clear
that $\wtvp{2}$ will collapse the image of a $3$-secant $2$-plane to 
a point, hence the exceptional locus must contain the transform of 
$Sec^2X$.  However by Theorem~\ref{embed}, we know that
the rational map $\varphi_2$ is an embedding off $Sec^1_3(SecX)$,
the trisecant variety to the secant variety.  This motivates the
following
\begin{lemma}\label{secantsareequal}
Let $X\subset \P^n$ be an irreducible variety.  Assume either of the following:
\begin{enumerate}
\item $Sec^kX$ is defined as a scheme by forms of degree $\leq 2k+1$.
\item $X$ is a smooth curve embedded by a line bundle of degree
at least $2g+2k+1$.
\end{enumerate}
Then $Sec^kX=Sec^1_{k+1}(Sec^{k-1}X)$ as schemes.
\end{lemma}

\begin{proof}

First, choose a $(k+1)$-secant $k$-plane $M$.  $M$ then intersects
$Sec^{k-1}X$ in a hypersurface of degree $k+1$, hence every line in
$M$ lies in $Sec^1_{k+1}(Sec^{k-1}X)$.  As $Sec^{k}X$ is
reduced and irreducible, 
$Sec^kX\subseteq Sec^1_{k+1}(Sec^{k-1}X)$ as schemes.

For the converse, assume the first condition is satisfied.
Choose a line $L$ that intersects 
$Sec^{k-1}X$ in a scheme
of length at least $k+1$.  It is easy to verify that $Sec^kX$ is
singular along $Sec^{k-1}X$, hence every form that vanishes on 
$Sec^kX$ must vanish $2k+2$ times on $L$.  By hypothesis, however,
$Sec^kX$ is scheme theoretically defined by forms of degree 
$\leq 2k+1$, hence each of these forms must vanish on $L$.

The sufficiency of the second condition follows from Thaddeus' 
construction and \cite[\S 2,(i)]{bertram3}.
\end{proof}

This implies that if $Sec^kX$ satisfies
$(K_{k+2})$ and if $Sec^kX=Sec^1_{k+1}(Sec^{k-1}X)$, then 
the map $\varphi_{k+1}:\P^{s_0}\dashrightarrow \P^{s_{k+1}}$ 
given by the forms defining $Sec^kX$ is an embedding off of $Sec^{k+1}X$.  
We use Theorem~\ref{embed} to understand the structure of these maps
via the following two lemmas:

\begin{lemma}\label{nolines}
If the embedding of a projective variety $X\subset\P^n$ is 
$(2k+4)$-very ample,  
then the intersection of two $(k+2)$-secant $(k+1)$-planes, if nonempty,
must lie in $Sec^kX$ (in fact, it must be an $\ell +1$ secant $\P^{\ell}$
for some $\ell \leq k$).  In particular, $Sec^{k+1}X$ has dimension
$(k+2)\dim X+k+1$.
\end{lemma}
\begin{proof}
The first statement is elementary: Assume two $(k+2)$-secant $(k+1)$-planes
intersect at a single point.  If the point is not on $X$, then there
are $2k+4$ points of $X$ that span a $(2k+2)$-plane, which is
impossible by hypothesis.  Hence the intersection lies
in $Sec^0X=X$.  A simple repetition of this argument for larger dimensional
intersections gives the desired result.
The statement of the dimension follows immediately; or see
\cite[11.24]{harris}.
\end{proof}
\begin{lemma}\label{possiblepreimages}
Let $X\subset \PP{0}$ be an irreducible variety whose embedding is 
$(2k+4)$-very ample.  Assume that $Sec^kX$ satisfies $(K_{k+2})$, and that 
$Sec^{k+1}X=Sec^1_{k+2}(Sec^kX)$ as schemes.  Let $\Gamma$ be the 
closure of the graph of
$\varphi_{k+1}$ with projection $\pi:\Gamma\rightarrow \PP{0}$.  If $a$ is
a point in the closure of the image of $\varphi_{k+1}$ and 
$F_a\subset \Gamma$ is the fiber over $a$ then
$\pi(F_a)$ is one of the following:
\begin{enumerate}
\item a reduced point in $\PP{0}\setminus Sec^{k+1}X$
\item a $(k+2)$-secant $(k+1)$-plane
\item {\bf contained} in a linear subspace of $Sec^kX$
\end{enumerate}
\end{lemma}

\begin{proof}

The first and third possibilities follow directly from 
Theorem~\ref{embed}.

For the second, note that a priori $\pi(F_a)$ could
be any linear space intersecting $Sec^kX$ in a hypersurface
of degree $k+2$.  However, Lemma~\ref{nolines} and the hypothesis
that $Sec^{k+1}X=Sec^1_{k+2}(Sec^kX)$ immediately imply that
any such linear space must be $k+1$ dimensional; hence a 
$(k+2)$-secant $(k+1)$-plane.
\end{proof}

With these results in hand we present the general construction.

Let $Y_0$ be an irreducible projective variety and suppose
$\beta_i:Y_0\dashrightarrow Y_i, 1\leq i\leq j$, is a 
collection of dominant, birational maps.
Define the {\bf dominating variety} of the collection, denoted
$\B_{(0,1,\ldots,j)}$, to be the 
closure of the graph of 
$$(\beta_1,\beta_2,\ldots,\beta_{j}):Y_0\dashrightarrow 
Y_1\times Y_2\times\cdots\times Y_j$$
Denote by $\B_{(a_1,a_2,\ldots,a_r)}$
the projection of $\B_{(0,1,\ldots,j)}$ to 
$Y_{a_1}\times Y_{a_2}\times\cdots\times Y_{a_r}$.
Note that $\B_{(a_1,a_2,\ldots,a_r)}$ is birationally isomorphic to
$\B_{(b_1,b_2,\ldots,b_k)}$ for all $0\leq a_r,b_k\leq j$.
Note further that if the $\beta_i$ are all morphisms
then $\B_{(0,1,\ldots,j)}\cong Y_0$, in other words only
rational maps contribute to the structure of the
dominating variety.

\begin{rem}{Definition/Notation}{k2l}
We say $X\subset \PP{0}$ satisfies condition
$(K_2^j)$ if $Sec^iX$ satisfies condition $(K_{2+i})$ for 
$0\leq i\leq j$; hence $X$ satisfies $(K_2^0)$ if and only
if $X$ satisfies $(K_2)$, $X$ satisfies $(K_2^1)$ if and only
if $X$ satisfies $(K_2)$ and $SecX$ satisfies $(K_3)$, etc.
{\nopagebreak \hfill $\Box$ \par \medskip}
\end{rem}
If $X\subset \PP{0}$ satisfies condition $(K_2^j)$, then each
rational map $\varphi_i:\PP{0}\dashrightarrow \PP{i}$ is
birational onto its image for $1\leq i\leq j+1$, and assuming 
the conclusion of Lemma~\ref{secantsareequal} each $\varphi_i$ is an
embedding off $Sec^iX$.  
Therefore $\B_{(i)}$ is the closure of the image of $\varphi_i$,
$\B_{(0,i)}$ is the closure of the graph of $\varphi_i$, and in 
the notation of Theorem~\ref{flip} $\B_{(0,1,2)}\cong\wt{M_2}$
and $\B_{(1,2)}\cong M_2$.  Note $\B_{(0)}=\PP{0}$.

\begin{lemma}\label{allblowups}
$\B_{(0,1,2,\ldots,i)}$ is the blow up of $\B_{(0,1,2,\ldots,i-1)}$
along the proper transform of $Sec^{i-1}X$, $1\leq i\leq j+1$.
\end{lemma}

\begin{proof}
This is immediate from the definition (or see \cite[3.1.1]{verm-thesis}).
\end{proof}

\begin{rem}{Remark}{bert-thadspaces}
The spaces constructed in \cite{bertram1} are of the type
$\B_{(0,1,2,\ldots,k)}$.  The spaces $\wt{M_k}$ and $M_k$ 
constructed in \cite{thad-flip}
are $\wt{M_k}\cong \B_{(k-2,k-1,k)}$ and 
$M_k\cong \B_{(k-1,k)}$.
{\nopagebreak \hfill $\Box$ \par \medskip}
\end{rem}

Our goal is to understand explicitly the geometry of this
web of varieties generalizing Theorem~\ref{flip}.
In the next section we describe in detail the structure
of the second flip.  As each subsequent flip requires the understanding
of $\H^kX$ for larger $k$, it is not clear that the process will
continue nicely beyond the second flip (at least for 
varieties of arbitrary dimension).

\section{Construction of the Second Flip}\label{secondflip}

Let $X\subset \PP{0}$ be a smooth, irreducible variety that 
satisfies $(K_2^1)$.
The diagram of varieties we study in this section is:
\begin{center}
{\begin{minipage}{1.5in}
\xymatrix{
 & & \B_{(0,1,2)} \ar[dl] \ar[dr]^{h_1} \ar[d] & & \B_{(1,2,3)} \ar@{..>}[dl] \ar@{..>}[dr]^{h_2} &   \\
 & \B_{(0,1)} \ar[dr]_{\wtvp{1}} \ar[dl] & \B_{(0,2)} \ar[dll] \ar[drr]  & \B_{(1,2)} \ar[dl] \ar[dr]^{\wtvp{2}} & & \B_{(2,3)} \ar@{..>}[dl]  \\
 \B_{(0)}  &   & \B_{(1)}  & & \B_{(2)} &   }
\end{minipage}}
\end{center}
where $\B_{(0,1,2)}$ is the dominating variety of the pair of birational
maps $\varphi_1:\PP{0}\dashrightarrow \PP{1}$ and 
$\varphi_2:\PP{0}\dashrightarrow \PP{2}$; and
where we have yet to construct the two rightmost varieties.  
We write $\picard \B_{(0,1)}=\picard \B_{(1,2)}=\Z H+\Z E$
and $\picard \B_{(0,1,2)}=\Z H+\Z E_1+\Z E_2$ (recall all three spaces 
are smooth by Theorem~\ref{flip}).

\begin{thm}\label{wtv2}
Let $X\subset \B_{(0)}=\PP{0}$ be a smooth, irreducible variety 
of dimension $r$ that satisfies $(K_2^1)$.
Assume that $X$ is embedded by a complete linear system $|L|$ 
and that the following conditions are satisfied:
\begin{enumerate}
\item $L$ is $(5+r)$-very ample
\item If $r\geq 2$, then for every point $p\in X$, $H^1(X,L\tensor\I_p^3)=0$
\item $Sec^2X=Sec^1_3(Sec^1X)$ as schemes
\item The projection of $X$ into $\P^{m}, m=s_0-1-r$, from
any embedded tangent space is such that the image is projectively normal
and satisfies $(K_2)$  
\end{enumerate}
Then the morphism $\wtvp{2}:\B_{(1,2)}\rightarrow \B_{(2)}$ induced by 
$\O_{\B_{(1,2)}}(3H-2E)$ 
is an embedding off the transform of $Sec^2X$, and the restriction of
$\wtvp{2}$ to the transform of $Sec^2X$ has fibers isomorphic to $\P^2$.
\end{thm}

As the proof of Theorem~\ref{wtv2} is somewhat involved, we break it into
several pieces.  We begin with a Lemma and a crucial observation, 
followed by the 
proof of the Theorem.  The observation invokes 
a technical lemma whose proof is postponed until the end.

\begin{rem}{Remark On the Hypotheses}{easilysatisfied}
Note that if $X$ is a smooth curve embedded by a line bundle of degree
at least $2g+5$, then conditions $1-4$ are automatically 
satisfied.  Conjecture~\ref{degreeboundsconj} would imply
condition $(K_2^1)$ holds also.  Furthermore, if $r=2$ 
and $H^1(X,L)=0$ then
condition $1$ implies condition $2$.

If $r\geq 2$, then the image of the projection from the space tangent
to $X$ at $p$ is $\blow{X}{p}\subset \P^{m}$.
Furthermore, by the discussion after Remark~\ref{5secant3plane} 
any such projection of $X$ will be generated as a scheme by quadrics when 
$SecX$ is defined by cubics, hence condition $4$ is not unreasonable.
{\nopagebreak \hfill $\Box$ \par \medskip}
\end{rem}

\begin{lemma}\label{projectionsarenotsobad}
With hypotheses as in Theorem~\ref{wtv2}, the image of the projection of 
$X$ into $\P^{m}$, $m=s_0-1-r$, is $\blow{X}{p}$, hence is smooth.  
Furthermore, it contains no lines and it contains no plane quadrics 
except for the exceptional divisor, which is the quadratic Veronese 
embedding of $\P^{r-1}$.  
\end{lemma}

\begin{proof}

If $r=1$ the statement is clear.  Otherwise, 
let $X'\subset \P^{m}$ denote the closure of the image of projection from the
embedded tangent space to $X$ at $p$.  As mentioned above, 
$X'\cong \blow{X}{p}$, hence is smooth.  
Let $E_p\subset X'$ denote the exceptional divisor.  The existence
of a line or plane quadric not contained in $E_p$ is
immediately seen to be impossible by the $(5+r)$-very
ampleness hypothesis.

As $\P^m=\P\Gamma(X',L\tensor \O(-2E_p))$ and as $E_p\cong \P^{r-1}$,
we have $L\tensor \O(-2E_p)|_{E_p}\cong \O_{\P^{r-1}}(2)$.
Condition $2$ implies this restriction is surjective on global sections.
\end{proof}

\begin{rem}{Observation}{crucialobs}
Let $\B_{(0,1,2)}\rightarrow \B_{(0)}$ be the projection and let $F_p$ be the
fiber over $p\in X$; hence $F_p$ is the blow up of
$\P^{m}$ along a copy of $\blow{X}{p}$.  
We again denote this variety by $X'\subset
\P^{m}$, and the embedding of $X'$ into $\P^{m}$ satisfies
$(K_2)$ by hypothesis.  The restriction of $\O_{\B_{(0,1,2)}}(3H-2E_1-E_2)$
to $F_p$ can thus be identified with $\O_{\blow{\P^{m}}{X'}}(2H'-E')$,
and, noting Lemma~\ref{projectionsarenotsobad}, it seems that 
Theorem~\ref{embed} could be applied.
Unfortunately, it is not clear that this restriction should be surjective
on global sections.  However, by Lemma~\ref{nearsurjrest} below, the 
image of the morphism on $F_p$ induced by the restriction of
global sections is isomorphic to the image of the morphism given by
the complete linear system $|\O_{\blow{\P^{m}}{X'}}(2H'-E')|$.
Hence by the fourth hypothesis and Lemma~\ref{projectionsarenotsobad}, 
the only collapsing that occurs in $F_p$ under the morphism
$\B_{(0,1,2)}\rightarrow \B_{(2)}$ is that of 
secant lines to $X'\subset \P^{m}$.  

Now, for some $p\in X$, suppose that a secant line $S$ in $F_p$ is collapsed
to a point by the projection $\B_{(0,1,2)}\rightarrow \B_{(2)}$.  
Then $S$ is the proper 
transform of a secant line to $X'\subset \P^m$,
but every such secant line is the intersection of $F_p$ with a
$3$-secant $\P^2$ through $p\in X$.  For example, 
if $S\subset F_p$ is the secant line through
$q,r\in X'$, $q,r\notin E_p$, then $S$ is the intersection of
$F_p$ with the proper transform of the plane spanned by $p,q,r$.
It should be noted that the two dimensional fiber associated to the
collapsing of a plane spanned by a quadric in the exceptional
divisor (Lemma~\ref{projectionsarenotsobad}) will take the place of
a $3$-secant $\P^2$ spanned by a non-curvilinear scheme 
contained in the tangent space at $p$.

Therefore, all the collapsing 
in the exceptional locus over a point $p\in X$ is 
associated to the collapsing of $3$-secant $2$-planes.  
{\nopagebreak \hfill $\Box$ \par \medskip}
\end{rem}

\begin{proof}(of Theorem~\ref{wtv2})
Let $a\in \B_{(2)}$ be a point in the image of $\wtvp{2}$.  
The fiber over $a$ is mapped isomorphically into $\B_{(1)}$ 
by the projection $\B_{(1,2)}\rightarrow \B_{(1)}$.  We are
therefore able to study $(\wtvp{2})^{-1}(a)$ by looking at the
fiber of the projection $\B_{(0,1,2)}\rightarrow \B_{(2)}$,
and projecting to $\B_{(0,1)}$ and to $\B_{(1)}$.

By applying Lemma~\ref{possiblepreimages} to the map 
$\B_{(0,2)}\rightarrow \B_{(2)}$, 
the projection to $\B_{(0,1)}$ is contained as a scheme in the 
{\em total} transform of one of the following (note the more refined 
division of possibilities): 
\begin{enumerate}
\item a point in $\PP{0}\setminus Sec^2X$
\item a $3$-secant $2$-plane to $X$ not contained in $SecX$
\item a linear subspace of $SecX$ not tangent to $X$
\item a linear subspace of $SecX$ tangent to $X$
\end{enumerate}

In the first case, there is nothing to show as the total transform of 
a point in $\PP{0}\setminus Sec^2X$ is simply a reduced point and the
map $\wtvp{1}$ to $\B_{(1)}$ is an embedding in a neighborhood of
this point.

If the projection is a $3$-secant
$2$-plane, then by Observation~\ref{crucialobs} the projection to $\B_{(0,1)}$ is 
a $3$-secant $2$-plane blown up at the three points of intersection, and
so the image in $\B_{(1)}$ is a $\P^2$ that has undergone a Cremona 
transformation.  

In the third case, Observation~\ref{crucialobs} shows that either 
the projection to $\B_{(0,1)}$ is the {\em proper} transform of a secant line
to $X$, or that the projection to $\B_{(0)}$ is a linear subspace of
$SecX$ that is not a secant line.  In the first case, every such
space is collapsed to a point by $\wtvp{1}$.  The second
implies $\wtvp{2}$ has a fiber of dimension $d$ that is contained in
$\P(\F)\subset \B_{(1,2)}$.  Because $E_2\rightarrow \P(\F)$ is a $\P^1$-bundle, this
implies the projection of the fiber to $\B_{(0)}$ is contained
in a linear subspace $M$ of $SecX$ of dimension $d+1$.  Furthermore,
the proper transform of $M$ is collapsed to a 
$d$ dimensional subspace of $\B_{(1)}$, in particular the general point
of $M$ lies on a secant line {\em contained in $M$} by Theorem~\ref{getabundle}.
Therefore $Y=M\cap X$ has $SecY=M$, hence $Sec^2Y=M$ but this is impossible
by Lemma~\ref{nolines} and the restriction that $M$ not be tangent
to $X$.  

In the final case, the proper transform in $\B_{(0,1)}$ of a linear 
space $M\cong \P^k$ tangent to $X$ at a point $p$ is $\blow{\P^k}{p}$.
Denote the exceptional $\P^{k-1}$ by $Q$; Lemma~\ref{projectionsarenotsobad} 
implies $Q\cong E_p$ is the quadratic Veronese embedding of 
$\P^{k-1}\subset \P(\Gamma(\blow{X}{p},L(-2E_p)))$.  
A simple dimension count 
shows that the restriction to $Q$ of the projective bundle $E_2\rightarrow \wts$ 
arising from the blow up of $\B_{(0,1)}$ along $\wts$ 
is precisely the restriction to $Q$ of the projective bundle 
arising from the {\em induced} blow up of $\P(\Gamma(\blow{X}{p},L(-2E_p)))$ along 
$\blow{X}{p}$; denote this variety $\P_Q$.  Furthermore, the transform of 
$\blow{\P^k}{p}$ in $\B_{(0,1,2)}$ is a $\P^1$-bundle over 
$\P_Q\subset\B_{(1,2)}$.  Now by Lemma~\ref{nearsurjrest}, every fiber
of $\wtvp{2}$ contained in $\P_Q\subset\B_{(1,2)}$ is either a point or
is isomorphic to a $\P^2$ spanned by a plane quadric in $Q$.
\end{proof}

\begin{rem}{Remark}{altprooffork=2}
For curves, parts $3$ and $4$ of the proof can also be concluded by showing
that any line contained in $SecX$ must be a secant or tangent line (this
is immediate from the $6$-very ample hypothesis).
{\nopagebreak \hfill $\Box$ \par \medskip}
\end{rem}

To complete the proof, we need Lemma~\ref{nearsurjrest}
which itself requires a general result:

\begin{lemma}\label{killideal}
Let $\pi:X\rightarrow Y$ be a flat morphism of smooth 
projective varieties.  Let $F=\pi^{-1}(p)$ be a smooth fiber and
let $L$ be a locally free sheaf on $X$.  If $R^i\pi_*L=0$ and
$H^i(F,L\tensor \O_F)=0$ for all $i>0$, 
then $R^i\pi_*(\I_F\tensor L)=0$ for all $i>0$.
\end{lemma}
\begin{proof}
The hypotheses easily give the vanishing
$R^i\pi_*(\I_F\tensor L)=0$ for all $i>1$.  For $i=1$, take the 
exact sequence on $Y$
$$0\rightarrow \pi_*(\I_F\tensor L)\rightarrow \pi_*L\rightarrow
\pi_*(\O_F\tensor L)\rightarrow R^1\pi_*(\I_F\tensor L)\rightarrow 0$$
Because $\pi_*(\O_F\tensor L)$ is supported at the point $p$, it
suffices to check that 
$H^1(F,\I_F\tensor \O_F\tensor L)=H^1(F,N^*_{F/X}\tensor L)=0$.
$\pi$ flat implies $N^*_{F/X}\cong \pi^*(N^*_{p/Y})$,
hence $N^*_{F/X}$ is trivial.  Now, $H^1(F,L\tensor \O_F)=0$ implies
$H^1(F,N^*_{F/X}\tensor L)=0$.
\end{proof}

\begin{lemma}\label{nearsurjrest}
Under the hypotheses of Theorem~\ref{wtv2}, The image of $F_p$
under the projection $\B_{(0,1,2)}\rightarrow \B_{(2)}$
is isomorphic to the image of $F_p$ under the morphism
induced by the complete linear system associated to
$\O_{F_p}(2H'-E')$.
\end{lemma}

\begin{proof}

{\bf Step 1:}  {\em If $a,b\in F_p$ are mapped to the same point under
the projection to $\B_{(2)}$, 
then $a$ and $b$ map to the same point under the projection to $\B_{(0,2)}$.}
This is clear from the construction of the maps in question as the projections
$\PP{0}\times\PP{1}\times\PP{2}\rightarrow \PP{2}$ and
$\PP{0}\times\PP{1}\times\PP{2}\rightarrow \PP{0}\times\PP{2}$ respectively.

{\bf Step 2:}  {\em Re-embed $\B_{(0,2)}\hookrightarrow \P^N\times\P^{s_2}$
via the map associated to $\O_{\PP{0}}(k)\bt\O_{\P^{s_2}}(1)$.}
This gives a map $\B_{(0,1,2)}\rightarrow \P^N\times\P^{s_2}$ induced
by a subspace of 
$H^0(\B_{(0,1,2)},\O_{\PP{0}}(k)\bt\O_{\P^{s_2}}(1)\tensor\O_{\B_{(0,1,2)}})$
where $\O_{\PP{0}}(k)\bt\O_{\P^{s_2}}(1)\tensor\O_{\B_{(0,1,2)}}\cong
\O_{\B_{(0,1,2)}}((k+3)H-2E_1-E_2)$.  
As $\B_{(0,2)}\hookrightarrow \P^N\times\P^{s_2}$ is an embedding,
the induced maps on $F_p$ have isomorphic images for all $k\geq 1$.  We have, therefore,
only to show $H^0(\PP{0}\times\P^{s_1}\times\P^{s_2},\O_{\PP{0}}(k)\bt\O_{\P^{s_2}}(1))$
surjects onto $H^0(F_p,\O_{F_p}(2H'-E'))$ for some $k$.

{\bf Step 3:} {\em The map $$H^0(\PP{0}\times\P^{s_1}\times\P^{s_2},\O_{\PP{0}}(k)\bt\O_{\P^{s_2}}(1))\rightarrow H^0(\B_{(0,1,2)},\O_{\B_{(0,1,2)}}((k+3)H-2E_1-E_2))$$ is surjective for
all $k\gg 0$.}  This follows directly from the fact that $SecX$ is scheme theoretically
defined by cubics and the construction of $\PP{2}$ as $\P(\Gamma(\B_{(0,1,2)},\O_{\B_{(0,1,2)}}(3H-2E_1-E_2)))$.

{\bf Step 4:} {\em The map $$H^0(\B_{(0,1,2)},\O_{\B_{(0,1,2)}}((k+3)H-2E_1-E_2))
\rightarrow H^0(E_1,\O_{E_1}((k+3)H-2E_1-E_2))$$ is surjective for
all $k\gg 0$.}  

We show $H^1(\B_{(0,1,2)},\O_{\B_{(0,1,2)}}((k+3)H-3E_1-E_2))=0$.  
Let $\rho:\B_{(0,1,2)}\rightarrow \B_{(0)}$ be the projection.  By the
projective normality assumption of Theorem~\ref{wtv2},
$R^i\rho_*\O_{E_1}((k+3)H-\ell E_1-E_2)=0$ for all $i,\ell>0$ since $E_1\rightarrow X$
is flat.  Ampleness of $\O_{\PP{0}}(H)$ implies
$H^1(E_1,\O_{E_1}(mH-\ell E_1-E_2))=0$ for all $m\geq m_0$, where $m_0$ may
depend on $\ell$.  From the exact sequence 
$$\ses{\O_{\B_{(0,1,2)}}(mH-(\ell+1)E_1-E_2)}{\O_{\B_{(0,1,2)}}(mH-\ell E_1-E_2)}{\O_{E_1}(mH-\ell E_1-E_2)}$$
a finite induction shows that if $H^1(\B_{(0,1,2)},\O_{\B_{(0,1,2)}}(mH-(\ell+1)E_1-E_2))=0$
for $m\gg 0$, some $\ell>1$ then $H^1(\B_{(0,1,2)},\O_{\B_{(0,1,2)}}((k+3)H-3E_1-E_2))=0$ for
all $k\gg 0$.

As $K_{\B_{(0,1,2)}}=\O_{\B_{(0,1,2)}}((-s_0-1)H+(s_0-r-1)E_1+(s_0-2r-2)E_2)$,
we have 
$$\O_{\B_{(0,1,2)}}(mH-(\ell+1)E_1-E_2-K)=\O_{\B_{(0,1,2)}}((m+s_0+1)H-(\ell+s_0-r)E_1-(s_0-2r-1)E_2)$$
As soon as $\ell\geq s_0-3r-2$, the right side is $\rho$-nef and, because $\rho$ is birational,
the restriction of the right side to the general fiber of $\rho$ is big.  Hence
by \cite[2.17.3]{kollar-pairs},
$R^i\rho_*\O_{\B_{(0,1,2)}}(mH-(\ell+1)E_1-E_2)=0$ for $i\geq 1$.  Again
by the ampleness of $\O_{\PP{0}}(H)$, we have $H^1(\B_{(0,1,2)},\O_{\B_{(0,1,2)}}(mH-(\ell+1)E_1-E_2))=0$
for $m\gg 0$, $\ell$ as above.

{\bf Step 5:} {\em The map $H^0(E_1,\O_{E_1}((k+3)H-2E_1-E_2))
\rightarrow H^0(F_p,\O_{F_p}(2H'-E'))$ is surjective for all $k\gg 0$.}
This is immediate by Lemma~\ref{killideal} and the projective
normality assumption of Theorem~\ref{wtv2}.
\end{proof}

As in Theorem~\ref{getabundle}, we show that the restriction of 
$\wtvp{2}$ to the
transform of $Sec^2X$ is a projective bundle over $\H^3X$.  
By a slight abuse of notation, write $\wt{Sec^2X}\subset \B_{(1,2)}$ for
the image of the proper transform of $Sec^2X$.  
Note the following:
\begin{lemma}\label{newrest}
Let $S_Z=(\wtvp{2})^{-1}(Z)\cong \P^2$ be a fiber over a point $Z\in\H^3X$.
Then $\O_{S_Z}(H)=\O_{\P^2}(2)$ and $\O_{S_Z}(E)=\O_{\P^2}(3)$.
\end{lemma}

\begin{proof}
This is immediate from the restrictions $\O_{S_Z}(2H-E)=\O_{\P^2}(1)$
and $\O_{S_Z}(3H-2E)=\O_{\P^2}$
\end{proof}

\begin{lemma}
There exists a morphism $\wt{Sec^2X}\rightarrow \G(2,s_0)$ whose
image is $\H^3X$.
\end{lemma}

\begin{proof}

A point $p\in\wt{Sec^2X}$ determines a unique
$2$-plane $S_Z$ in $\wt{Sec^2X}$ by Theorem~\ref{wtv2}.   
For every such $p$, the homomorphism
$H^0(\B_{(1,2)},\O_{\B_{(1,2)}}(H))\rightarrow H^0(\B_{(1,2)},\O_{S_Z}(H))$
has rank $3$, hence gives a point in $\G(2,s_0)$.  The image of the
associated morphism clearly coincides with the natural embedding 
of $\H^3X$ into $\G(2,s_0)$ described in \cite{cat-gott}.
\end{proof}

As in \cite[3.5]{verm-flip1}, there is a morphism $\H^3X
\rightarrow \B_{(2)}$ so that the composition factors $\wtvp{2}:
\wt{Sec^2X}\rightarrow \B_{(2)} $.  This is constructed by associating to 
every $Z\in\H^3X$ the rank $1$ homomorphism:
$$H^0(\B_{(1,2)},\O_{\B_{(1,2)}}(3H-2E))\rightarrow 
H^0(\B_{(1,2)},\O_{S_Z}(3H-2E))$$
where $S_Z$ is the $\P^2$ in $\B_{(1,2)}$ associated to $Z$. 

Exactly as in Theorem~\ref{getabundle}, this allows the identification of
$\wt{Sec^2X}$ with a $\P^2$-bundle over $\H^3X$.  Specifically,
$\E_2=(\wtvp{2})_*(\O_{\wt{Sec^2X}}(2H-E))$ is a rank $3$ vector 
bundle on $\H^3X$ and:
\begin{prop}
With notation as above, $\wtvp{2}:\wt{Sec^2X}\rightarrow \H^3X$ is the
$\P^2$-bundle $\P_{\H^3X}(\E_2)\rightarrow \H^3X$.
{\nopagebreak \hfill $\Box$ \par \medskip}
\end{prop}

We wish to show further that blowing up $Sec^2X$ along $X$ and then 
along $SecX$ resolves the singularities of $Sec^2X$.  By 
Theorem~\ref{flip}, $h_1:\B_{(0,1,2)}\rightarrow \B_{(1,2)}$ is 
the blow up of $\B_{(1,2)}$ along $\P(\F)$, hence it
suffices to show $\P(\F)\cap \P(\E_2)$ is a smooth subvariety
of $\P(\E_2)$.

\begin{prop}\label{intersectbundles}
$\D=\P(\F)\cap \P(\E_2)$ is the nested Hilbert scheme
$Z_{2,3}(X)\subset \H^2X\times \H^3X$, hence is smooth.  Therefore 
$\blow{\blow{Sec^2X}{X}}{\wts}\subset \B_{(0,1,2)}$ is smooth and
$Sec^2X\subset \PP{0}$ is normal.
\end{prop}

\begin{proof}
Let $\U_i\subset X\times \H^iX$
denote the universal subscheme.
We have morphisms $\wtvp{2}:\D\rightarrow \H^3X$ and 
$\wtvn{1}:\D\rightarrow \H^2X$, and it is routine to check that
$(id_X\times\wtvn{1})^{-1}(\U_2)\subset (id_X\times\wtvp{2})^{-1}(\U_3)$.
Hence (Cf. \cite[\S 1.2]{lehn}) $\wtvn{1}\times\wtvp{2}$ maps $\D$ to 
the nested Hilbert scheme $Z_{2,3}(X)\subset \H^2X\times \H^3X$, where
closed points of $Z_{2,3}(X)$ correspond to pairs of subschemes $(\alpha,\beta)$
with $\alpha\subset \beta$.  Furthermore, via the description
of the structure of the map $\wtvp{2}$, it is clear that the
morphism of $\H^3X$-schemes $\D\rightarrow Z_{2,3}(X)$ is finite
and birational.  It is shown in \cite[0.2.1]{cheah} that
$Z_{2,3}(X)$ is smooth, hence this is an isomorphism.
\end{proof}

Let $\B_{(1,2,3)}$ be the blow up of $\B_{(1,2)}$ along
$\P(\E_2)$; note $\B_{(1,2,3)}$ is smooth.  
To construct $\B_{(2,3)}$, we first construct
the exceptional locus as a projective bundle over $\H^3X$.
Write $\picard\B_{(1,2,3)}=\Z H+\Z E_1+\Z E_3$.

\begin{lemma}\label{buildF2}
Let $p_3:E_3\rightarrow \H^3X$ be the composition $E_3\rightarrow
\P(\E_2) \rightarrow \H^3X$.  Then $\F_2=(p_3)_*\O_{E_3}(4H-3E_1-E_3)$ is
locally free of rank $s_0-3r-2=codim~(\wt{Sec^2X},\B_{(1,2)})$.
\end{lemma}

\begin{proof}

Each fiber $F_x$ of $p_3$ is isomorphic to $\P^2\times \P^t$, 
$t+1=codim~(\wt{Sec^2X},\B_{(1,2)})$.  Furthermore  
$H^0(F_x,\O_{F_x}(4H-3E_1-E_3))=H^0(\P^t,\O_{\P_t}(1))$
follows easily from Lemma~\ref{newrest}.
\end{proof}

There is a map $E_3\rightarrow \P(\F_2)$ given by the surjection
$$p_3^*\F_2\rightarrow \O_{E_3}(4H-3E_1-E_3)\rightarrow 0$$
hence a diagram of exceptional loci:
\begin{center}
{\begin{minipage}{1.5in}
\xymatrix{
 & E_3 \ar[dl] \ar[dr] \ar[dd]^{p_3} & \\
 \P(\E_2) \ar[dr] & & \P(\F_2) \ar[dl] \\
 & \H^3X &  }
\end{minipage}}
\end{center}
It is important to note that 
$$\P(\F_2)\cong \P (p_{3_*}\O_{E_3}(4H-3E_1-E_3+t(3H-2E_1)))$$
for all $t\geq 0$ as the direct image on the right will differ from $\F_2$
by a line bundle.  Hence for all $t\geq 0$ the same morphism 
$E_3\rightarrow \P(\F_2)$ is induced by the surjection
$$p_3^*p_{3_*}\O_{E_3}(4H-3E_1-E_3+t(3H-2E_1))\rightarrow 
\O_{E_3}(4H-3E_1-E_3+t(3H-2E_1))$$

One can now repeat almost verbatim \cite[4.7-4.10]{verm-flip1}
to construct the second flip; i.e. the space $\B_{(2,3)}$.
Recall the following:
\begin{prop}\cite[4.5]{verm-flip1}\label{free}
Let $\L$ be an invertible sheaf on a complete variety $X$, and 
let $\mathscr{B}$ be any locally free sheaf.
Assume that the map $\lambda:X\rightarrow Y$ 
induced by $|\mathscr{L}|$ is a birational morphism and that $\lambda$ 
is an isomorphism in a neighborhood of $p\in X$.
Then for all $n$ sufficiently large, the map
$$H^0(X,\mathscr{B}\tensor \mathscr{L}^n)\rightarrow
H^0(X,\mathscr{B}\tensor \mathscr{L}^n\tensor \O_p)$$
is surjective.
{\nopagebreak \hfill $\Box$ \par \medskip}
\end{prop}

Taking $\mathscr{B}=\O_{\B_{(1,2,3)}}(4H-3E_1-E_3)$ and 
$\L=\O_{\B_{(1,2,3)}}(3H-2E_1)$, the map
induced by the linear system associated to
$$\O_{\B_{(1,2,3)}}((4H-3E_1-E_3)+(k-2)(3H-2E_1))= 
\O_{\B_{(1,2,3)}}((3k-2)H-(2k-1)E_1-E_3) $$
is base point free off $E_3$ for $k\gg 3$.  To show this gives a morphism,
one shows the restriction of above linear system to the divisor $E_3$ induces
a surjection on global sections, hence restricts to the map
$E_3\rightarrow\P(\F_2)$ above.  For this, define
$\L_{\rho}=\O((3\rho -2)H-(2\rho -1)E_1-E_3)$ and write 
$$\O_{\B_{(1,2,3)}}((3k-2)H-(2k-1)E_1-2E_3)\tensor K_{\B_{(1,2,3)}}^{-1}=
\L_{\alpha}^{s_0-3r-1}\tensor \mathcal{A} $$
where $\alpha =\frac{2k+s_0-r-1}{2s_0-6r-2}$ and 
$\mathcal{A}=\O_{\B_{(1,2,3)}}\left( \left( \frac{3s_0-9r-4}{2}\right) H-(s_0-3r-2)E_1\right) $.
By the above discussion, $\L_{\alpha}^{s_0-3r-1}$ is nef for $k\gg 0$
and it is routine to verify that $\mathcal{A}$ is a big and nef
$\Q$-divisor; hence $H^1(\B_{(1,2,3)},\O((3k-2)H-(2k-1)E_1-2E_3))=0$.

The variety $\B_{(2,3)}$ is defined to be the image of this morphism.
This gives:

\begin{prop}\label{topofflip2}
With hypotheses as in Theorem~\ref{wtv2} and for $k$ sufficiently large, 
the morphism $h_2:\B_{(1,2,3)}\rightarrow\B_{(2,3)}$ induced by the linear
system $|\mathscr{L}_k|$ is an embedding off of $E_3$ and
the restriction of $h_2$ to $E_3$ is the morphism
$E_3\rightarrow \P(\F_2)$ described above. 
{\nopagebreak \hfill $\Box$ \par \medskip}
\end{prop}

\begin{rem}{Remark}{bestforh3}
The best (smallest) possible value for $k$ is $k=3$.  This will
be the case if $Sec^3X\subset \PP{0}$ is scheme theoretically
cut out by quartics.
{\nopagebreak \hfill $\Box$ \par \medskip}
\end{rem}

\begin{lemma}
$\B_{(2,3)}$ is smooth.
\end{lemma}

\begin{proof}
Because $\B_{(2,3)}$ is the image of a smooth variety with reduced, connected fibers
it is normal (Cf. \cite[3.2.5]{verm-thesis}).  Let $Z\cong \P^2$ be a fiber
of $h_2$ over a point $p\in \P(\F_2)$.  $Z\times\{ p\} $ is a fiber of a 
$\P^2\times \P^t$ bundle over $\H^3X$, hence the normal bundle sequence
becomes:
$$\ses{\bigoplus_1^{s_0-3} \O_{\P^2}}{N_{Z/\B_{(1,2,3)}}}{\O_{\P^2}(-1)}$$
This sequence splits, and allowing the elementary calculations
$H^1(Z,S^rN_{Z/\B_{(1,2,3)}})=0$ and 
$H^0(Z,S^rN_{Z/\B_{(1,2,3)}})=S^rH^0(Z,N_{Z/\B_{(1,2,3)}})$ for all $r\geq 1$, 
$\B_{(2,3)}$ is smooth by a natural extension of the smoothness portion of Castelnuovo's
contractibility criterion for surfaces given in \cite[2.4]{aw}.
\end{proof}

Letting $\P(\F_0)=\P_X(N^*_{X/\PP{0}})=E_1$ and $\P(\E_0)=X$,
the analogue of Theorem~\ref{flip} is:
\begin{thm}\label{flip2}
Let $X\subset \B_{(0)}=\PP{0}$ be a smooth, irreducible variety 
of dimension $r$ that satisfies $(K_2^1)$, with $s_0\geq 3r+4$.
Assume that $X$ is embedded by a complete linear system $|L|$ 
and that the following conditions are satisfied:
\begin{enumerate}
\item $L$ is $(5+r)$-very ample and $Sec^2X=Sec^1_3(Sec^1X)$ as schemes
\item The projection of $X$ into $\P^{m}, m=s_0-1-r$, from
any embedded tangent space is such that the image is projectively normal
and satisfies $(K_2)$  
\item If $r\geq 2$, then for every point $p\in X$, $H^1(X,L\tensor\I_p^3)=0$
\end{enumerate}
Then there is a pair of flips as pictured below with:
\begin{enumerate}
\item $\B_{(i,i+1)}$ and $\B_{(i,i+1,i+2)}$ smooth
\item $\B_{(i,i+1)}\setminus \P(\E_{i+1})\cong \B_{(i+1,i+2)}\setminus \P(\F_{i+1})$;
as $s_0\geq 3r+4$, $\picard\B_{(0,1)}\cong \picard\B_{(i+1,i+2)}$
\item $\P\E_i=\P \wtvp{i_*}\O_{\wt{Sec^iX}}(iH-(i-1)E)$ and $\P\F_i=\P \wtvn{i_*}\O_{\wt{Sec^iX}}((i+2)H-(i+1)E)$
\item $h_i$ is the blow up of $\B_{(i,i+1)}$ along $\P(\F_i)$
\item $\B_{(i,i+1,i+2)}\rightarrow \B_{(i,i+1)}$ is the blow up along $\P(\E_{i+1})$
\item $\wtvn{i}$, induced by $\O_{\B_{(i,i+1)}}((i+1)H-iE)$, is an embedding off
of $\P(\F_i)$, and the restriction of $\wtvn{i}$ is the projection $\P(\F_i)\rightarrow \H^{i+1}X$ 
\item $\wtvp{i}$, induced by $\O_{\B_{(i-1,i)}}((i+1)H-iE)$,
is an embedding off of $\P(\E_i)$, and the restriction of
$\wtvp{i}$ is the projection $\P(\E_i)\rightarrow \H^{i+1}X$
\item $\P(\F_i)\cap\P(E_{i+1})\subset \B_{(i,i+1)}$ is isomorphic to
the nested Hilbert scheme $Z_{i+1,i+2}\subset \H^iX\times \H^{i+1}X$,
hence is smooth.
\end{enumerate}
\end{thm}

\begin{center}
{\begin{minipage}{1.5in}
\xymatrix{
 & & & \B_{(0,1,2)} \ar[dl] \ar[dr]^{h_1} & & \B_{(1,2,3)} \ar[dl] \ar[dr]^{h_2} &   \\
\B_{(0)}\ar@{=}[dr]_{\wtvp{0}} & & \B_{(0,1)} \ar[dr]_{\wtvp{1}} \ar[dl]^{\wtvn{0}} & & \B_{(1,2)} \ar[dl]^{\wtvn{1}} \ar[dr]_{\wtvp{2}} & & \B_{(2,3)} \ar[dl]^{\wtvn{2}}  \\
 & \B_{(0)}  &   & \B_{(1)}  & & \B_{(2)} &   }
\end{minipage}}
\end{center}

\begin{center}
{\begin{minipage}{1.5in}
\diagram
 & E_{j+1} \dlto \drto^{h_j} & \\
  \P(\E_j) \drto_{\wtvp{j}} & & \P(\F_j) \dlto^{\wtvn{j}} \\
 & \H^{j+1}X & 
\enddiagram
\end{minipage}}
\end{center}
{\nopagebreak \hfill $\Box$ \par \medskip}

\end{document}